\documentclass{statsoc}

\usepackage{graphicx}
\usepackage{natbib}

\title{Voting power and Qualified Majority Voting with a ``no vote'' option}

\author{Martin Kurth}
\address{School of Mathematical Sciences, University of Nottingham,
         UK.}
\email{martin.kurth@nottingham.ac.uk}
\coaddress{Martin Kurth, School of Mathematical Sciences, University of Nottingham, University Park, Nottingham, NG7 2RD, UK.}

\begin{document}

\begin{abstract}
In recent years, enlargement of the European Union has led to increased interest in the
allocation of voting weights to member states with hugely differing population numbers.
While the eventually agreed voting scheme lacks any strict mathematical basis, the Polish government
suggested a voting scheme based on the Penrose definition of voting power, leading to an
allocation of voting weights proportional to the square root of the population (the ``Jagiellonian
Compromise''). The Penrose definition of voting power is derived from the citizens' freedom
to vote either ``yes'' or ``no''. This paper defines a corresponding voting power based on
``yes'', ``no'' and ``abstain'' options, and it is found that this definition also leads to
a square root law, and to the same optimal vote allocation as the Penrose scheme.

\keywords{Penrose voting; Qualified Majority Voting; Square root voting; Voting power; Voting systems} 
\end{abstract}

\section{Introduction}

Following the failure of the draft EU constitution \citep{DraftConstitution} in referenda in France and the
Netherlands, and the subsequent negotiations on a Reform Treaty \citep{ReformTreaty}, the voting arrangements
in the Council of the EU have received a significant amount of public attention, not least
through the Polish proposal of a voting system that gives every member state a voting
weight proportional to the square root of its population \citep{euobserver}. The idea of
this voting scheme is to give every EU citizen the same influence on decisions in the
Council, based on an analysis of their voting power. The concept of voting power used
in this context was first introduced by \citet{penrose1}, and adapted to
the EU framework by \citet{slomczynski1}, \citet{slomczynski3}.
In this scheme for Qualified Majority Voting, the threshold for a motion to pass in the Council of the EU is
then set according to an optimality condition, see \citet{slomczynski1}, \citet{slomczynski3},
\citet{Kurth1} for details.

The Penrose definition of voting power assumes that all citizens have the freedom
to vote either ``yes'' or ``no'' in elections or referenda. However, this definition
of voting power ignores the freedom not to vote at all, which is a freedom very
frequently used by the citizens of EU member states, as voter turnout in elections
is generally much less than 100\%.

The purpose of this paper is to provide an analogous definition of voting power which
includes this ``no vote'' option.
In section~\ref{section21}, a brief introduction to the Penrose voting power scheme is
given, while in section~\ref{section22} it is generalised to the case where citizens
have the possibility to abstain from voting. A short conclusion summarises the result.

\section{Voting power}

Wherever constituencies with different numbers of voters elect representatives for
a council where the representatives for each constituency have to cast a joint vote,
the question of allocation of voting weights to the constituencies arises. This is,
for example, the case in the European Union, where population numbers vary between
about 400,000 (Malta) and 82,300,000 (Germany), see \citet{EuroStat}. While it may be intuitive to allocate
voting weights proportional to the population numbers, there is actually no sound
mathematical foundation for such a scheme. What is needed is a rigorous definition
of citizens' voting power.

\subsection{Penrose voting power}
\label{section21}

Penrose defined voting power as the probability that the vote of a single voter is
decisive, provided all other voters vote randomly, in an election where voters
choose between ``yes'' and ``no''.

More precisely, let us assume there are $N+1$ voters in the country, and a motion is
successful if it receives more ``yes'' than ``no'' votes. Whether or not the vote
of a single voter is decisive depends on how the remaining $N$ voters have voted.
For a single voter to be decisive, the motion must be successful if they vote ``yes'', and
fail if they vote ``no''. If $N$ is even, this is the case when the $N$ votes are
evenly split between ``yes'' and ``no'', if $N$ is odd, the individual voter will
be decisive if there is one ``yes'' vote more than there are ``no'' votes.

Penrose's fundamental assumption is that these $N$ voters vote randomly, ie all
$2^{N}$ possible voting outcomes are equally likely.

The number of voting outcomes
where the remaining voter is decisive is given by the binomial coefficients
\begin{equation}
\left(\begin{array}{c}N\\N/2\end{array}\right)\quad (N\,\,{\rm even}),\qquad
\left(\begin{array}{c}N\\(N-1)/2\end{array}\right)\quad (N\,\,{\rm odd}).
\end{equation}

This means that the probability for the remaining voter to be decisive, ie their
\emph{voting power} $P_N$, is
\begin{equation}
P_N=\frac{1}{2^{N}}\left(\begin{array}{c}N\\ \lfloor N/2\rfloor\end{array}\right),
\end{equation}
where $\lfloor N/2\rfloor$ denotes the floor function, ie the largest integer not
exceeding $N/2$.
Using the Stirling approximation for the factorials we obtain
\begin{equation}
P_N\approx\sqrt{\frac{2}{\pi}}\frac{1}{\sqrt{N}}\approx\sqrt{\frac{2}{\pi}}\frac{1}{\sqrt{N+1}}\qquad(N\gg 1),
\end{equation}
ie the voting power decreases like one divided by the square root of the population.

\subsection{Voting power with a ``no vote'' option}
\label{section22}

Let us now consider a voting scheme where all voters have the choice between ``yes'', ``no'' and
``abstain''. Let us again assume that there are $N+1$ voters, and a motion is successful if
it receives more ``yes'' than ``no'' votes. As in the Penrose  case, we define
the voting power of a single voter as the probability for this voter to be decisive, provided
the remaining $N$ voters vote randomly.

There are $3^{N}$ possible voting outcomes, all of them equally likely according to the random
voting assumption. For any integer $0\leq K\leq N$, there are
\[
 \left(\begin{array}{c}N\\K\end{array}\right)
\]
ways to choose $K$ voters participating, ie voting either ``yes'' or ``no''. If $K$ is
even, the single remaining voter will be decisive if the numbers of ``yes'' and ``no'' votes
are equal, either by voting ``yes'' and rendering the motion successful, or by voting ``no'' or
abstaining, in which case the motion will fail. If $K$ is odd, the single voter will be
decisive if there is one ``yes'' vote more than there are ``no'' votes; in this case a
``yes'' vote or an abstention will cause the motion to be successful, while a ``no'' vote
will lead to its rejection. This means that for any choice of $K$ voters who have voted
either ``yes'' or ``no'', there are
\[
 \left(\begin{array}{c}K\\ \lfloor K/2\rfloor\end{array}\right)
\]
voting outcomes rendering the single remaining voter decisive. The probability for this
one voter to be decisive, ie is their voting power, is thus
\begin{equation}
 \label{exactvotingpower}
P_N = \frac{1}{3^{N}}\sum_{K=0}^N\left(\begin{array}{c}N\\K\end{array}\right)
 \left(\begin{array}{c}K\\ \lfloor K/2\rfloor\end{array}\right).
\end{equation}

There does not seem to be any simple expression for $P_N$. What we are interested in here is its behaviour
for large $N$, as in the Penrose case.
For $N>0$, the sum can be written as
\begin{equation}
  \sum_{K=0}^N\left(\begin{array}{c}N\\K\end{array}\right) \left(\begin{array}{c}K\\ \lfloor K/2\rfloor\end{array}\right)
  = \left(\begin{array}{c}N\\0\end{array}\right)_2+\left(\begin{array}{c}N\\1\end{array}\right)_2
\end{equation}
\citep{callan1},
where $(\cdot)_2$ denotes the trinomial coefficient defined by
\begin{equation}
(x^2+x+1)^n=\sum_{k=-n}^n\left(\begin{array}{c}n\\k\end{array}\right)_2 x^{n+k}.
\end{equation}
The trinomial coefficients have asymptotic expansions
\begin{equation}
 \left(\begin{array}{c}N\\0\end{array}\right)_2 =
 \frac{\sqrt{3}}{2}\left(\begin{array}{c}2N\\N\end{array}\right)\left(\frac{3}{4}\right)^N
 \left(1+O\left(\frac{1}{N}\right)\right)
\end{equation}
and
\begin{equation}
 \left(\begin{array}{c}N\\1\end{array}\right)_2 =
 \frac{\sqrt{3}}{6}\left(\begin{array}{c}2(N+1)\\N+1\end{array}\right)\left(\frac{3}{4}\right)^{N+1}
 \left(1+O\left(\frac{1}{N}\right)\right)
\end{equation}
\citep{Merlini1}.
Using Stirling's formula for the binomial coefficients, we get
\begin{eqnarray}
 P_N & = & \frac{1}{2}\sqrt{\frac{3}{\pi}}\left(\frac{1}{\sqrt{N}}+\frac{1}{\sqrt{N+1}}\right)\left(1+O\left(\frac{1}{N}\right)\right) \\
\label{squarerootapprox}
 & \approx & \sqrt{\frac{3}{\pi}}\frac{1}{\sqrt{N+1}},
\end{eqnarray}
for large $N$,
which means that even in the case of a ``no vote'' option, the voting power of each citizen decreases
with the square root of the population number, as in the case of pure yes-no voting.

\begin{figure}
\input{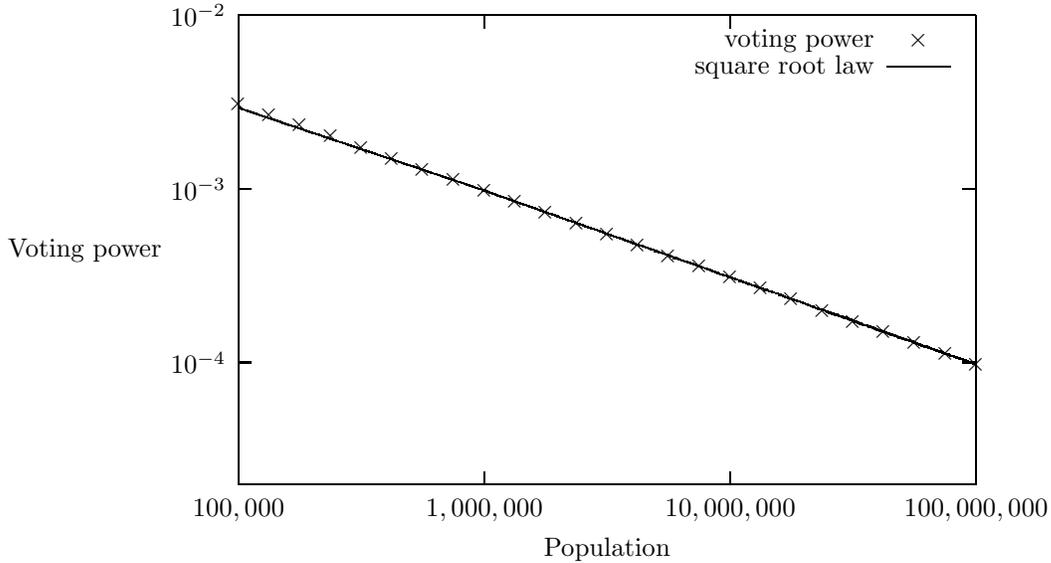}
\caption{Voting power for a range of population numbers. The crosses represent the voting power calculated from
equation~\ref{exactvotingpower}, the line is the square root approximation from equation~\ref{squarerootapprox}.}
\label{figure1}
\end{figure}

For some reasonable population numbers (ranging from $100,000$ to $100$ million), figure~\ref{figure1} shows the voting
power calculated from equation~\ref{exactvotingpower}, and compares it with the square root approximation
from equation~\ref{squarerootapprox}. It can be seen that for these population numbers, the error from truncating
the asymptotic expansions is negligible, and the square root approximation provides a good measure for
voting power.

\section{Conclusion}

In the previous section, it has been shown that in the case of voters having three voting options,
``yes'', ``no'' and ``abstain'', the voting power defined as the probability for a single voter to
be decisive is proportional to one divided by the square root of the population, in the limit of
large population numbers. This analysis is based on the random voting assumption, ie the assumption
that all possible voting outcomes are equally likely. While this assumption has been criticised in
the literature, eg in \citet{Garrett1}, \citet{Albert1}, \citet{bafumi1}, as not being realistic, 
other authors, eg \citet{slomczynski1}, \citet{slomczynski3}, \citet{Leech2}, \citet{Shapley1}, \citet{banzhaf2},
\citet{Coleman1}, point out that this is an \emph{a priori} voting power based on
voters' possible options rather than their actual decisions in previous elections.

The square root behaviour of the voting power is similar to voting power in pure yes-no votes, the only
difference being the constant factor. As the basis for voting weight allocation in the Jagiellonian Compromise
scheme for the EU is relative rather than absolute voting power, it is not affected by the constant
factor. This means that even the inclusion of a ``no vote'' option in the analysis of citizens' voting
power leads to the allocation of exactly the same voting weights under the Jagiellonian Compromise scheme.

\bibliography{paper}

\end{document}